\documentclass[11pt]{article}
\usepackage{a4}
\usepackage{amsmath}
\usepackage{amssymb}
\usepackage{amsthm}
\usepackage{amsfonts}
\usepackage{graphicx}
\usepackage[numbers, sort&compress]{natbib}
\usepackage{cite}
\usepackage{enumerate}

\usepackage{cleveref}

\newtheorem{theorem}{Theorem}
\newtheorem{lemma}{Lemma}

\newtheorem{definition}{Definition}

\newtheorem{proposition}{Proposition}

\begin{document}
\def\R{{\mathbb R}}
\def\Z{{\mathbb Z}}
\def\N{{\mathbb N}}
\def\S{{\mathbb{S}^{d-1}}}
\def\wh{\widehat}
\def\C{{\mathbb C}}
\def\lV{\left\Vert}
\def\rV{\right\Vert}
\def\lv{\left\vert}
\def\rv{\right\vert}
\title{Strict positive definiteness of convolutional and axially symmetric kernels on $d$-dimensional spheres}

\author{ Martin Buhmann \thanks{Email: martin.buhmann@math.uni-giessen.de} and Janin J\"ager \thanks{ Corresponding author. Email:
janin.jaeger@math.uni-giessen.de}\\
Justus-Liebig University, Lehrstuhl Numerische Mathematik, \\ Heinrich-Buff Ring 44, 35392 Giessen, Germany   }
\date{\today}

\maketitle

\begin{abstract}{The paper introduces new sufficient conditions of strict positive definiteness for kernels on $d$-dimensional spheres which are not radially symmetric but possess specific coefficient structures. The results use the series expansion of the kernel in spherical harmonics. The kernels either have a convolutional form or are axially symmetric with respect to one axis. The given results on convolutional kernels generalise the result derived by Chen et al.\ \citep{Chen2003} for radial kernels.\\}{Keywords: strictly positive definite kernels,\ covariance functions,\ sphere }
\end{abstract}

\section{Introduction}

Spherical approximation is a topic of immense interest and the use of positive definite spherical basis functions has recently been discussed in a tremendous number of publications \citep{Gneiting2013,Hubbert2015, Chen2003, Nie2018}. Most of the known results are for isotropic
positive definite kernels, which are kernels that only depend on the
geodesic distance of their arguments.  Isotropic kernels are the generalisation of the radial basis function method to the sphere and therefore referred to in approximation theory as spherical
radial basis functions (see \citep{Hubbert2015} and references therein). They
are also of importance in statistics where they occur as correlation
functions of homogeneous random fields on spheres \citep{Lang2015,
  Ma2015}. 
  
Recently  the use of axially symmetric kernels on the sphere
was suggested and applied for the approximation of global data in
\citep{Bissiri2019, Emery2019}.  Conditions that ensure strict positive definiteness of such kernels on the 2-sphere  were proven in \citep{Bissiri2020}. These results are in this paper generalised to $d$-dimensional spheres and additional necessary conditions are proven. 

The second type of non-radial kernels we study, are convolutional kernels which have for example been described in \citep{Schaback2001},\citep{Jaeger2021}. 

We will briefly summarize necessary definitions in the first section
and then, starting from a general form of the kernel, state conditions which ensure certain properties of the kernel, like axial symmetry, convolutional form and invariance under parity.  In the third section, we give sufficient conditions for axially symmetric kernels  to be positive definite and strictly positive definite and add some necessary conditions.  In the fourth part of this paper 
we generalise the result of Chen et al.\  \citep{Chen2003} on radial kernels to convolutional kernels deriving sufficient conditions on strict positive definiteness and finally study the special case of the circle in the last section.

\subsection{Problem description 
and background}

We focus on interpolation problems on the $d$-sphere 
$$\S=\left \lbrace \xi \in \R^d \vert\  \xi_1^2+\xi_2^2+\cdots +\xi_d^2=1 \right \rbrace,\quad d\geq2,$$
 where a finite set of distinct data sites $\Xi\subset \mathbb{S}^d$ and  values $f\left(\xi\right) \in \C$, $\xi \in \Xi$, of a possibly elsewhere unknown function $f$ on the sphere are given.

The approximant is formed as a linear combination of kernels \[K:\S \times \S \rightarrow \C.\] Taking the form 
\begin{equation}\label{eq:Interpolant}  s_f\left(\zeta\right)=\sum_{\xi\in \Xi} c_{\xi} K\left(\zeta, \xi\right), \qquad  \zeta \in \S,
\end{equation}
the  problem of finding such an approximant $s_f$ satisfying
 \begin{equation}
 s_f\left(\xi\right)=f\left(\xi\right), \qquad \forall \xi \in \Xi,
 \end{equation}
 is uniquely solvable under certain conditions on $K$. We assume all the kernels to be Hermitian, meaning they satisfy $K\left(\xi,\zeta\right)=\overline{K\left(\zeta, \xi\right)}$, so that the positive definiteness of the kernel will ensure the solvability of the interpolation problem for arbitrary data sets.
\begin{definition}\label{DF:SPD}
A Hermitian kernel $K: \Omega \times \Omega \rightarrow \C$ is called {\/\rm
positive definite on} $\Omega$ if  the matrix  $K_{\Xi}=\left\lbrace
K\left(\xi,\zeta\right)\right\rbrace_{\xi,\zeta \in \Xi}$ is positive semi-definite on $\C^{\vert \Xi\vert}$ for arbitrary finite sets of distinct points $\Xi\subset \Omega$.

The kernel is {\/\rm strictly positive definite} if $K_{\Xi}$ is a positive definite matrix on $\C^{\vert \Xi\vert}$  for  arbitrary finite sets of distinct points $\Xi$.
\end{definition} 

We assume that $K$ is continuous in both arguments, it is thereby square-integrable and can be represented as 
\begin{equation}\label{eq:GenKer} K\left(\xi, \zeta\right)=\sum_{j,j'=0}^{\infty}\sum_{k=1}^{N_{j,d}} \sum_{k'=1}^{N_{j',d}} a_{j,j',k,k'}Y_j^k\left(\xi\right)\overline{Y_{j'}^{k'}\left(\zeta\right)},\qquad \forall \xi,\zeta \in \S ,
 \end{equation}
where the $Y_j^k$ form an orthonormal basis of the eigenfunctions of the Laplace-Beltrami operator on the sphere and the corresponding eigenvalues increase with $j$. The $a_{j,j',k,k'}$ are complex numbers. To be precise the eigenfunctions are the solutions of the eigenfunction problem
$$\label{eq:EigenProb}\lambda f +\triangle f=0, $$
where $\triangle$ is the Laplace-Beltrami operator on the sphere.
For example on the $2$-sphere such a basis can take the form
\begin{equation}\label{eq:2spherhar}
Y_j^k\left(\theta,\varphi\right)= \frac{1}{\sqrt{2\pi}} \sqrt{\frac{2j+1}{2} \frac{\left(j-k\right)!}{\left(j+k\right)!}} P_{j}^k\left(\cos\left(\theta\right)\right)e^{ik \varphi},
\end{equation}
 and $P_{j}^k$ are the associated Legendre polynomials (see for example 8.1.2, \citep{Abramowitz1972}).

The eigenfunctions corresponding to the eigenvalues $\lambda_{j}=j\left(j+d-1\right)$ are spherical harmonics and the number of eigenfunctions corresponding to the eigenvalue $\lambda_j$ is denoted by $N_{j,d}$. The numbers are given by $N_{0,d}=1$,
\[N_{j,d}=\frac{\left(2j+d-2\right)\left(j+d-3\right)!}{j!\left(d-2\right)!}.\] 

We will denote the space of all spherical harmonics corresponding to the eigenvalue $\lambda_j$ by $H_j:=\operatorname{span} \{ Y_{j}^k,\ k=1,\ldots,N_{j,d}\}$.

Additionally we will use the following estimate which holds for any orthonormal basis  and is a direct consequence of the Poisson summation formula (\citep{Atkinson2012}, Equation (2.35))
\begin{equation}\label{eq:abssphH}
\vert Y_j^{k}\left(\xi\right)\vert \leq  \sqrt{\frac{N_{j,d}}{\sigma_{d-1}}},\quad  \forall \xi \in \S,
\end{equation}
where $\sigma_{d-1}$ is the surface area of $\S$.

It is well known and used in the characterisation above that every function in $L^2\left(\S\right)$ can be represented as a spherical harmonic expansion of the form 
\begin{equation}\label{eq:Serrep}
g\left(\xi\right)= \sum_{j=0}^{\infty}\sum_{k=1}^{N_{j,d}}\hat{g}_{j,k}Y_{j}^k\left(\xi\right),\qquad \text{ with } \hat{g}_{j,k}=\int_{\S}g\left(\xi\right)\overline{Y_{j}^k\left(\xi\right)}\,d \sigma\left( \xi\right),
\end{equation}
where $d\sigma$ is the surface area measure on $\S$.
For this expansion the Parseval equation for spherical harmonics  holds (\citep{Atkinson2012}, (2.143)): 
$$\Vert g\Vert_{L^2\left(\S\right)}^2=\sum_{j=0}^{\infty}\sum_{k=1}^{N_{j,d}} \vert \hat{g}_{j,k}\vert ^2.$$



\section{Kernels with special coefficient structure}

We now fix a basis for our space of spherical harmonics such that the kernels which are still Hermitian and continuous are represented as in \eqref{eq:GenKer}.
Imposing certain conditions on the structure of the coefficients allows us to focus on kernels with specific properties. 

The case studied most often is assuming that the kernel is isotropic \citep{Schoenberg1942b}, which means it only depends on the distance of its two arguments and not on their position on the sphere. 
The coefficients of an isotropic positive definite kernel given in the form  \eqref{eq:GenKer} satisfy
$$a_{j,j',k,k'}=c_j \delta_{j,'j} \delta_{k,k'}\geq 0$$ 
as stated by Schoenberg  \citep{Schoenberg1942b}. A characterisation of the strictly positive definite kernels of this form was presented by Chen et al.\ in \citep{Chen2003}. 

Recently, a sufficient condition  for axially symmetric kernels on the two-sphere  to be
strictly positive definite was presented in \citep{Bissiri2020}. The
axially symmetric kernels do depend on the difference in longitude of
the two inputs  $\xi,\zeta$ and their individual values of latitude.
The coefficients of these kernels
satisfy $$a_{j,j',k,k'}=c_k\left(j,j'\right)\delta_{k,k'},$$ 
where the basis of the spherical harmonics  on the 2-sphere is chosen as \eqref{eq:2spherhar}. The condition was first stated in \citep{Jones1963}.

For the study of properties of the kernel we need the kernel of the form \eqref{eq:GenKer} to be uniquely recoverable from its Fourier expansion on $\S \times \S$.
\begin{lemma}\label{Lemma:UniqueRep}
Let $K$ be a kernel of the form \eqref{eq:GenKer}. If 
\begin{equation}\label{eq:Critabssum} \sum _{j,j'=0}^{\infty}\sum_{k,k'=1}^{N_{j,d},N_{j',d}} \vert a_{j,j',k,k'}\vert \sqrt{\frac{N_{j,d}N_{j',d}}{\sigma_{d-1}^2}} <\infty, 
\end{equation}
then, for each combination of  $j,j'\in \N,$ $1\leq k \leq N_{j,d}, 1\leq  k'\leq N_{j',d}$ the coefficients $a_{j,j',k,k'}$ can be uniquely determined through
\begin{equation}\label{eq:FouriercoefKer}
a_{j,j',k,k'}=\int_{\S}\int_{\S}K\left(\xi,\zeta\right)Y_j^k\left(\xi\right) \overline{Y}_{j'}^{k'}\left(\zeta\right)\, d\sigma\left(\xi\right)\,d\sigma\left(\zeta\right).
\end{equation}
\end{lemma}
\begin{proof} The result follows by inserting the kernel representation \eqref{eq:GenKer} into  \eqref{eq:FouriercoefKer} and then using the estimate \eqref{eq:abssphH} to exchange the order of summation and integration since it implies absolute summability.
\end{proof}

First we note 
\begin{lemma} Let $K$ be a kernel of the form \eqref{eq:GenKer} which satisfies \eqref{eq:Critabssum}. It is Hermitian if and only if 
\begin{equation}\label{eq:CoefConv} a_{j,j',k,k'}=\overline{a_{j',j,k',k}}
\end{equation}
for all possible choices $j,j'\in \N$ and $1\leq k \leq N_{j,d}, 1\leq  k'\leq N_{j',d}.$
\end{lemma}
\begin{proof}
Assuming $K$ is Hermitian implies that 
$$\overline{a_{j',j,k',k}}=\int_{\S}\int_{\S}\overline{K\left(\xi,\zeta\right)}\overline{Y}_{j'}^{k'}\left(\xi\right) {Y}_{j}^{k}\left(\zeta\right)\, d\sigma\left(\xi\right)\,d\sigma\left(\zeta\right)=a_{j,j',k,k'}$$
according to \eqref{eq:FouriercoefKer}. 
The other direction follows from representation \eqref{eq:GenKer}.
\end{proof}
Besides axially symmetric kernels we are interested in kernels which are referred to as convolutional kernels. For a fixed basis of spherical harmonics these are expressed as kernels with an eigenvalue block structure.

\begin{lemma} A Hermitian kernel  $K\in L^2\left(\S\times\S\right)$ of the form \eqref{eq:GenKer} satisfies 
\[\int_{\S}K\left(\xi,\zeta\right)\overline{Y\left(\xi\right)}d\sigma(\xi) \in H_j, \quad \forall\ Y\in H_j,\]
if and only if \[a_{j,j',k,k'}=\delta_{j,j'}d_j\left(k,k'\right).\]
\end{lemma}
\begin{proof}
It is sufficient to prove the implication for the basis $Y_{j}^k$, $k=1,\ldots,N_{j,d}$, because of the linearity of the integral.
For these
\[\int_{\S}K\left(\xi,\zeta\right)\overline{Y_{j}^k\left(\xi\right)}d\sigma (\xi)=\sum_{j'=0}^{\infty} \sum_{k'=1}^{N_{j',d}}a_{j,j',k,k'}Y_{j'}^{k'}\left(\zeta\right),\]
which is an element of $H_j$ if and only if $a_{j,j',k,k'}=\delta_{j,j'}d_j\left(k,k'\right)$.
\end{proof}

Kernels of this form are invariant under parity, meaning $K\left(\xi,\zeta\right)=K\left(-\xi,-\zeta\right)$, and the special structure allows to determine easily if an interpolant derived using such a kernel is included in certain  Sobolev spaces, as for example studied in \citep{Narcowich2002}.

The kernels have been mostly discussed without the selection of a fixed basis of the spherical harmonics of a certain order. In this case one can choose the basis of $H_j$, denoted by $ \tilde{Y}_j^k\left(\xi\right)$ in a way that 
 \begin{equation}\label{eq:Kernform3} K\left(\xi, \zeta\right)=\sum_{j=0}^{\infty}\sum_{k=1}^{N_{j,d}}  d_{j,k} \tilde{Y}_j^k\left(\xi\right)\overline{\tilde{Y}_{j}^{k}\left(\zeta\right)},\qquad \forall \xi,\zeta \in \S.
 \end{equation}
 Under  the above conditions and for $f\in L^2\left(\S\right)$ the  convolution operator with kernel $K$ is
 $$Tf\left(\xi\right)=\int_{\S}K\left(\xi,\zeta\right)f\left(\zeta\right)d\sigma\left(\zeta\right),$$
 and the Fourier coefficients of $Tf$ as in \eqref{eq:Serrep} when computed with respect to the basis $\tilde{Y}_j^k$ satisfy  $\widehat{\tilde{Tf}}_{j,k}=d_{j,k} \hat{\tilde{f}}_{j,k}$.
 We nevertheless include the slightly more complex block structure, since we believe it to be helpful when we want to construct new kernels using our results. For example we may be using the fixed basis of the spherical harmonics on the 2-sphere as mentioned. 

For completeness we include the conditions on the coefficients which ensure invariance under parity. One can easily verify that all restrictions of shift-invariant, Hermitian kernels in $\R^d$ to the surface of the sphere are invariant under parity.
\begin{lemma}[Parity-invariance]
A kernel of the form \eqref{eq:GenKer} with property \eqref{eq:Critabssum}  satisfies,  $$K\left(\xi,\zeta\right)=K\left(-\xi,-\zeta\right),\quad \forall \xi,\zeta\in\S$$ if and only if
\begin{equation}
a_{j,j',k,k'}=0, \quad \forall j+j'\neq 0\mod 2.
\end{equation}
\end{lemma}
\begin{proof}
The results follows immediately using the property of homogenity of the spherical harmonics:
\begin{align*}
K\left(-\xi,-\zeta\right)=\sum_{j=0}^{\infty}\sum_{k=1}^{N_{j,d}} \sum_{j'=0}^{\infty}\sum_{k'=1}^{N_{j,d}} a_{j,j',k,k'}\left(-1\right)^{j+j'}\overline{Y_{j'}^{k'}}\left(\zeta\right)Y_{j}^{k}\left(\xi\right), \quad \forall \xi,\zeta \in \S.
\end{align*}
This is equal to $K\left(\xi,\zeta\right)$ if and only if all coefficients are equal, meaning $ a_{j,j',k'k'}=0$,  for all pairs of $j,j'$ with odd sum.
\end{proof}

The previous conditions are independent of the choice of orthogonal basis $Y_j^k$. For axial symmetry we will have to fix a certain basis, for the existing results for $d=3$ this was \eqref{eq:2spherhar}
where axial symmetry with respect to $\varphi$ led to the coefficients satisfying $a_{j,j',k,k'}=\delta_{k,k'}c_{j,j'}$.
On general spheres with $d\geq 2$ we will focus on axial symmetry with respect to just one of the $d-1$ axes and we assume that the sphere has been rotated such that this is the axis which is represented by the first coordinate in the polar coordinate form of $\xi \in \S$, $\xi=\left(\theta_1,\ldots,\theta_{d-1}\right)^T$. The spherical harmonics of degree $\ell_{d-1}$ can explicitly be given by
\begin{equation}\label{eq:spherharmaxial}
Y_{\ell_1,\ldots,\ell_{d-1}}\left(\theta_1,\ldots,\theta_{d-1}\right)=\frac{1}{\sqrt{2\pi}}e^{i\ell_1\theta_1} \prod_{j=2}^{d-1}{}_{j}\tilde{P}_{\ell_j}^{\ell_{j-1}}\left(\theta_{j}\right),\end{equation}
where $\ell_1,\ldots,\ell_{d-1}$ are integers satisfying
\[\ell_{d-1}\geq \cdots \geq \vert \ell_1\vert\]
and 
$${}_{j}\tilde{P}_{L}^{\ell}(\theta)={}_{j}c_L^{\ell}\left(\sin\left(\theta\right)\right)^{-\left(2-j\right)/2}P_{L+\left(j-2\right)/2}^{-\left(\ell+\left(j-2\right)/2\right)}\left(\cos\left(\theta\right)\right),$$
where $P_{\nu}^{\mu}$ are the associated Legendre functions and
$${}_{j}c_L^{\ell}:= \left( \frac{2L+j-1}{2}\frac{\left(L+\ell+j-2\right)!}{\left(L-\ell\right)!}\right)^{1/2}.$$
The formula is taken from \citep{Higuchi1987}, Equation (2.5).
To be able to handle the formula better  we define the last part of the spherical harmonics as $\underline{\ell}:=\left(\ell_2,\ldots,\ell_{d-2}\right)$, $\theta':= \left(\theta_2,\ldots,\theta_{d-1}\right)$ and
\[ p_ {\underline{\ell},\ell_{d-1}}\left(\theta'\right)= \frac{1}{\sqrt{2\pi}}\prod_{j=2}^{d-1}{}_{j}\tilde{P}_{\ell_j}^{\ell_{j-1}}\left(\theta_{j}\right).\]
Further we set $\Lambda_{\ell_1,\ell_{d-1}}:=\left\lbrace \left(\ell_2,\ldots, \ell_{d-2}\right)\in \N^{d-3} \vert\, \vert \ell_1\vert \leq \ell_2 \leq \cdots \leq \ell_{d-2}\leq \ell_{d-1} \right \rbrace$.
\begin{theorem}[Axial symmetry with respect to one axis]\label{THM:CoefStrucAX}
A kernel of the form 
\begin{align*}K\left(\xi,\zeta\right)=&\sum_{\ell_{d-1},\ell_{d-1}'=0}^{\infty}\ \sum_{\ell_1=-\ell_{d-1}}^{\ell_{d-1}} \sum_{\ell_1'=-\ell_{d-1}'}^{\ell_{d-1}'}\ \sum_{\underline{\ell}\in \Lambda_{\ell_1,\ell_{d-1}}} \ \sum_{\underline{\ell'}\in \Lambda_{\ell_1',\ell_{d-1}'}} a_{\ell_1,\underline{\ell},\ell_{d-1},\ell_1',\underline{\ell'},\ell_{d-1}'} \times \\
&\times Y_{\ell_1,\underline{\ell},\ell_{d-1}}\left(\xi\right)\overline{Y_{\ell_1',\underline{\ell'},\ell'_{d-1}}}\left(\zeta\right),\quad \xi,\zeta \in \S,  \end{align*}
 which satisfies \eqref{eq:Critabssum} is axially symmetric with respect to the $\theta_1$ axis if and only if 
\begin{equation}\label{eq:AxialCoef} a_{\ell_1,\ldots,\ell_{d-1},\ell_1',\ldots,\ell_{d-1}'}= \delta_{\ell_1,\ell'_1}c_ {\ell_1}(\ell_2,\ldots\ell_{d-1},\ell_2',\ldots,\ell_{d-1}').
\end{equation}
\end{theorem}
\begin{proof} First we simplify the expression for the kernel as follows, with $\xi=\left(\theta_1,\ldots,\theta_{d-1}\right)$ and $\zeta=\left(\upsilon_1,\ldots,\upsilon_{d-1}\right)$ given in polar coordinates:
\begin{align*} K\left(\xi,\zeta\right)=&\sum_{\ell_{d-1},\ell_{d-1}'=0}^{\infty}\ \sum_{\ell_1=-\ell_{d-1}}^{\ell_{d-1}}\ \sum_{\ell_1'=-\ell_{d-1}'}^{\ell_{d-1}'}\ \sum_{\underline{\ell}\in \Lambda_{\ell_1,j},\underline{\ell}'\in \Lambda_{\ell_1',\ell'_{d-1}} } a_{\ell_1,\underline{\ell},\ell_{d-1},\ell_1',\underline{\ell}',\ell'_{d-1}}\\
&\times e^{-i\ell_1\theta_1+i\ell_1'\upsilon_1} p_ {\underline{\ell},\ell_{d-1}}\left(\theta'\right) \overline{p_ {\underline{\ell'},\ell'_{d-1}}\left(\upsilon'\right)} \\
=&\sum_{\ell_{d-1},\ell_{d-1}'=0}^{\infty}\ \sum_{\ell_1=-\ell_{d-1}}^{\ell_{d-1}}\ \sum_{\ell_1'=-\ell'_{d-1}}^{\ell'_{d-1}} a_{\ell_1,\underline{\ell},\ell_{d-1},\ell_1',\underline{\ell}',\ell'_{d-1}}e^{-i\ell_1\theta_1+i\ell_1'\upsilon_1} \\
& \times \sum_{\underline{\ell}\in \Lambda_{\ell_1,\ell_{d-1}},\underline{\ell}'\in \Lambda_{\ell_1',\ell'_{d-1}} }  p_ {\underline{\ell},\ell_{d-1}}\left(\theta'\right) \overline{p_ {\underline{\ell'},\ell_{d-1}}\left(\upsilon'\right)}.
\end{align*}
The rotation in the $\theta_1$-axis by the angle $\alpha$ will be denoted as $R_{\alpha}$. Then 
\begin{align*}
K\left(R_{\alpha}\xi,R_{\alpha}\zeta\right)&=\sum_{\ell_{d-1},\ell'_{d-1}=0}^{\infty}\ \sum_{\ell_1=-\ell_{d-1}}^{\ell_{d-1}}\sum_{\ell_1'=-\ell'_{d-1}}^{\ell'_{d-1}}a_{\ell_1,\underline{\ell},\ell_{d-1},\ell_1',\underline{\ell}',\ell'_{d-1}}e^{-i\ell_1\left(\theta_1+\alpha\right)+i\ell_1'\left(\vartheta_1+\alpha\right)} \\
&\times \sum_{\underline{\ell}\in \Lambda_{\ell_1,\ell_{d-1}},\underline{\ell}'\in \Lambda_{\ell'_1,\ell'_{d-1}} }  p_ {\underline{\ell},\ell_{d-1}}\left(\theta'\right) \overline{p_ {\underline{\ell'},\ell_{d-1}}\left(\upsilon'\right)}\\
&=\sum_{\ell_{d-1},\ell'_{d-1}=0}^{\infty}\ \sum_{\ell_1=-\ell_{d-1}}^{\ell_{d-1}}\ \sum_{\ell_1'=-\ell'_{d-1}}^{\ell'_{d-1}}\  \sum_{\underline{\ell}\in \Lambda_{\ell_1,\ell_{d-1}},\underline{\ell}'\in \Lambda_{\ell_1',\ell'_{d-1}} } e^{-i\alpha\left(\ell_1-\ell'_1\right)} \\
&\times a_{\ell_1,\underline{\ell},\ell_{d-1},\ell_1',\underline{\ell'},\ell_{d-1}'}Y_{\ell_1,\ldots,\ell_{d-1}}\left(\xi\right)\overline{Y_{\ell_1',\ldots,\ell'_{d-1}}}\left(\zeta\right).
\end{align*}
We see that the above condition is sufficient for axial symmetry. It is also necessary since according to \Cref{Lemma:UniqueRep} equality is only possible if the coefficients are all equal and this can only hold for all $\alpha>0$ if 
$a_{\ell_1,\underline{\ell},\ell_{d-1},\ell_1',\underline{\ell'},\ell_{d-1}'}= \delta_{\ell_1,\ell'_1}c_ {\ell_1}(\ell_2,\ldots\ell_{d-1},\ell_2',\ldots,\ell_{d-1}').$
\end{proof}
As expected, the axially symmetric kernels could be written as a function depending on the values $\theta'$,$\upsilon'$ and the longitudinal difference $\theta_1-\upsilon_1$.
In \citep{Emery2019} additionally the properties of being longitudinal independent and longitudinal-reversible are introduced.
\begin{definition}
\begin{itemize}
\item An axially symmetric kernel is called longitudinal reversible if 
\[K\left(\left(\theta_1,\theta'\right),\left(\upsilon_1,\upsilon'\right)\right)=K\left(\left(\upsilon_1,\theta'\right),\left(\theta_1,\upsilon'\right)\right), \quad \forall \theta_1, \upsilon_1\in[0,2\pi]. \]
\item An axially symmetric kernel is called longitudinal-independent if 
\[K\left(\left(\theta_1,\theta'\right),\left(\varphi_1,\varphi'\right)\right)=K\left(\left(\theta,\theta'\right),\left(\varphi,\varphi'\right)\right), \quad \forall \theta_1,\theta, \varphi_1,\varphi \in[0,2\pi]. \]
\end{itemize}
\end{definition}For the next lemma and further use in the next section we introduce a further simplified notation.
First we define for an axially symmetric, Hermitian kernel
$$c_{\ell_1}: \Lambda_{\vert \ell_1\vert }\times \Lambda_{\vert \ell_1\vert } \rightarrow \C,$$ where 
$\Lambda_{\vert \ell_1\vert }=\left \lbrace \left(\ell_2,\ldots,\ell_{d-2},\ell_{d-1}\right)\in \N^{d-2}\vert\ \vert \ell_1\vert \leq \ell_2 \leq \cdots \leq \ell_{d-2}\leq \ell_{d-1}\right\rbrace$ 
and
\begin{equation}\label{eq:CoeffuctAx} c_{\ell_1}\left((\ell_2,\ldots,\ell_{d-2},\ell_{d-1}),(\ell_2',\ldots,\ell'_{d-2},{\ell'}_{d-1}\right)),
\end{equation}
where the coefficients are as in \eqref{eq:AxialCoef}. We deduce from \Cref{THM:CoefStrucAX} that every Hermitian, continuous and axially symmetric kernel which satisfies \eqref{eq:Critabssum}
can therefore be written as 
\begin{equation}\label{eq:Formaxker2}
 K\left(\xi, \zeta\right)=\sum_{\ell_1=-\infty}^{\infty}\sum_{\ell,\ell' \in \Lambda_{\vert \ell_1\vert }}c_{\ell_1}\left(\ell,\ell'\right) Y_{\ell_1,\ell}\left(\xi\right) Y_{\ell_1,\ell'}\left(\zeta\right),
\end{equation}
where $Y_{\ell_1,\ell}\left(\xi\right)$, with $\ell=\left(\ell_2,\ldots,\ell_{d-1}\right)\in\Lambda_{\vert \ell_1\vert}$ is the abbreviation for $Y_{\ell_1,\ell_2,\ldots,\ell_{d-1}}\left(\xi\right)$. Additionally, we find that these spherical harmonics can be represented for $\xi=\left(\theta_1,\ldots,\theta_{d-1}\right)$ in polar coordinates
\begin{equation}
Y_{\ell_1,\ell}\left(\xi\right)=\begin{cases} e^{i\ell_1\theta_1} g_{\ell}\left(\theta'\right),& \ell_1\geq0 \\
(-1)^{\ell_1}e^{i\ell_1\theta_1} g_{\ell}\left(\theta'\right),& \ell_1<0\end{cases}
\end{equation}
where $\theta'=\left(\theta_2,\ldots, \theta_{d-1}\right)$, $g_{\ell}$ is real valued, and the representation follows from \eqref{eq:spherharmaxial}.

We briefly state sufficient conditions of the expansion coefficients for longitudinal reversibility and independence. 
\begin{lemma}
\begin{enumerate}
\item An axially symmetric kernel satisfying \eqref{eq:Critabssum} is longitudinal reversible if 
\[ c_{\ell_1}\left(\ell,\ell'\right)=c_{-\ell_1}\left(\ell,\ell'\right),\quad \forall \ell_1\in \Z, \ \ell,\ell'\in \Lambda_{\vert \ell_1\vert}.\]
\item An axially symmetric kernel satisfying \eqref{eq:Critabssum} is longitudinal independent if $c_{\ell_1}\left(\ell,\ell'\right)=0$ for all $\ell_1\neq0.$
\item An axially symmetric kernel satisfying \eqref{eq:Critabssum} with the above expansion is real valued if and only if 
\begin{equation}\label{EQ:RealAxial}
c_{\ell_1}\left(\ell,\ell'\right)=\overline{c_{-\ell_1}\left(\ell,\ell'\right)},\quad \forall \ell_1\in \Z,\ \ell,\ell'\in\Lambda_{\vert \ell_1\vert}.
\end{equation}
\end{enumerate}
\end{lemma}
\begin{proof}

1. Using the above representation of the spherical harmonics for $\zeta=\left(\upsilon_1,\ldots,\upsilon_{d-1}\right)$,  in polar coordinates, with $\upsilon'=\left(\upsilon_{2},\ldots,\upsilon_{d-1}\right)$, we find  
\begin{align*}
 K\left(\xi, \zeta\right)=&\sum_{\ell_1=-\infty}^{\infty}\sum_{\ell,\ell' \in \Lambda_{\vert \ell_1\vert }}c_{\ell_1}\left(\ell,\ell'\right)e^{i\ell_1\theta_1-i\ell_1\upsilon_1} g_{\ell}\left(\theta'\right)g_{\ell'}\left(\upsilon'\right)\\
 =&\sum_{\ell_1=1}^{\infty}\sum_{\ell,\ell' \in \Lambda_{\vert \ell_1\vert }}\left( e^{i\ell_1(\theta_1-\upsilon_1)}c_{\ell_1}\left(\ell,\ell'\right)+e^{-i\ell_1(\theta_1-\upsilon'_1)}c_{-\ell_1}\left(\ell,\ell'\right)\right) g_{\ell}\left(\theta'\right)g_{\ell'}\left(\upsilon'\right)\\
 &+  \sum_{\ell,\ell' \in\Lambda_{0}}c_0\left(\ell,\ell'\right) g_{\ell}\left(\theta'\right)g_{\ell'}\left(\upsilon'\right).
\end{align*}
Setting as assumed $c_{\ell_1}\left(\ell,\ell'\right)=c_{-\ell_1}\left(\ell,\ell'\right)$, the above simplifies to  
\begin{align*}
 K\left(\xi, \zeta\right)=2&\sum_{\ell_1=1}^{\infty}\sum_{\ell,\ell' \in \Lambda_{\vert \ell_1\vert }}c_{\ell_1}\left(\ell,\ell'\right)\cos\left(\ell_1\left(\theta_1-\upsilon_1\right)\right)g_{\ell}\left(\theta'\right)g_{\ell'}\left(\upsilon'\right)\\
 &+  \sum_{\ell,\ell' \in\Lambda_{0}}c_0\left(\ell,\ell'\right) g_{\ell}\left(\theta'\right)g_{\ell'}\left(\upsilon'\right),
\end{align*}
where the values of $\theta_1$ and $\upsilon_1$ can be exchanged without changing the value of the kernel because of the symmetry of the cosine.

2. The sufficiency of the second part follows directly from the representation in the last equation, which shows that if all $c_{\ell_1}$ are constant zero for $\ell_1\neq0$ the kernel is independent of the values $\theta_1$, $\upsilon_1$.

3. We use the same representation above for the kernels $K\left(\xi,\zeta\right)$ as well as for its complex conjugate:

\begin{align*} \overline{K\left(\xi, \zeta\right)}=&\sum_{\ell_1=-\infty}^{\infty}e^{-i\ell_1\left(\theta_1-\upsilon_1\right)} \sum_{\ell,\ell' \in \Lambda_{\vert \ell_1\vert }}\overline{c_{\ell_1}\left(\ell,\ell'\right)} g_{\ell}\left(\theta'\right)g_{\ell'}\left(\upsilon'\right)
\end{align*}
and deduce that because of the uniqueness of the expansion coefficients the two are equal if and only if $\overline{c_{\ell_1}\left(\ell,\ell'\right)}=c_{-\ell_1}\left(\ell,\ell'\right)$ for all $\ell_1\in\Z$ and $\ell,\ell'\in \Lambda_{\vert \ell_1\vert}$.
\end{proof}In general a real basis of the spherical harmonics can be chosen if only real kernels should be studied. Using a real basis of the spherical harmonics will result in slightly different conditions on the coefficients. 

We finally summarise the connections between the introduced properties, which follows directly from the definitions.
\begin{lemma}
We assume the kernels in this proposition to be Hermitian and to satisfy \eqref{eq:Critabssum}.
\begin{enumerate}
\item In the case $d=2$ all axially symmetric kernels have convolutional form.
\item An axially symmetric kernel given in the form \cref{eq:Formaxker2} has a convolutional form if and only if $$c_{\ell_1}(\underline{\ell},\ell_{d-1},\underline{\ell'},\ell'_{d-1})=0$$ for all $\ell_{d-1}\neq \ell'_{d-1}\in\N$, .
\end{enumerate}
\end{lemma}

\section{(Strict) positive definiteness of axially symmetric kernels}
In \citep{Bissiri2020} strictly positive definite axially symmetric kernels on the $2$-sphere were described and sufficient conditions for strict  positive definiteness stated. We add sufficient conditions for $d$-dimensional spheres and prove some additional necessary conditions. In this section we focus on kernels on spheres with $d\geq3$, the special case of the circle is discussed in Section 5.

From now on we will use property \eqref{eq:Critabssum} also for axially symmetric and convolutional kernels. In this case the coefficients $a_{j,j',k,k'}$ in \eqref{eq:Critabssum} are derived using \eqref{eq:CoeffuctAx} and \eqref{eq:AxialCoef} for axially symmetric kernels and \eqref{eq:CoefConv} in the convolutional case. 
\begin{theorem}\label{Le:MapAxiaPD}
Let $K$ be an Hermitian, axially symmetric kernel of the form \eqref{eq:Formaxker2}, with $d\geq3$, satisfying \eqref{eq:Critabssum} for absolute summability.
The kernel is positive definite if the mapping $c_{\ell_1}: \Lambda_{\vert \ell_1\vert }\times \Lambda_{\vert \ell_1\vert} \rightarrow \C$ is positive definite for all $\ell_1 \in \Z$.
\end{theorem}
\begin{proof}We rewrite the quadratic form 
\begin{align*} \sum_{\xi,\zeta \in \Xi}\lambda_{\xi} \overline{\lambda_{\zeta}} K\left(\xi, \zeta\right)=&\sum_{\ell_1=-\infty}^{\infty}\sum_{\ell,{\ell'}\in \Lambda_{\vert \ell_1\vert}}c_{\ell_1}\left(\ell,\ell'\right) \sum_{\xi \in \Xi}\lambda_{\xi}Y_{\ell_1,\ell}\left(\xi\right) \overline{\sum_{\zeta \in \Xi}\lambda_{\zeta}Y_{\ell_1,\ell'}\left(\zeta\right)} .
\end{align*}
Now we define $y^{\ell_1,\underline{\ell},\ell_{d-1}}_{\Xi}:=\sum_{\xi \in \Xi}\lambda_{\xi} Y_{\ell_1,\underline{\ell},\ell_{d-1}}\left(\xi\right) \in\C$ to find 
\begin{align*}
 \sum_{\xi,\zeta \in \Xi}\lambda_{\xi} \overline{\lambda_{\zeta}} K\left(\xi, \zeta\right)=& \sum_{\ell_1=-\infty}^{\infty}\underset{k \rightarrow \infty}{\lim} \sum_{\ell_{d-1},{\ell'}_{d-1}=\vert \ell_1\vert }^{k}\ \sum_{\underline{\ell}\in \Lambda_{\ell_1,\ell_{d-1}}}\ \sum_{\underline{\ell'}\in \Lambda_{\ell_1,{\ell'}_{d-1}}}\\
&\times  y^{\ell_1,\underline{\ell},\ell_{d-1}}_{\Xi} c_{\ell_1}\left(\left(\underline{\ell},\ell_{d-1}\right),\left(\underline{\ell'},{\ell'}_{d-1}\right)\right) \overline{y^{\ell_1,\underline{\ell}',{\ell'}_{d-1}}_{\Xi}},\\
\end{align*}
where the absolute summability is necessary to ensure equality after reordering.
The positive definiteness of the mapping $c_{\ell_1}$ implies that
$$\sum_{\ell_{d-1},\ell'_{d-1}=\vert \ell_1\vert }^{k}\ \sum_{\underline{\ell}\in \Lambda_{\ell_1,\ell_{d-1}}} \ \sum_{\underline{\ell'}\in \Lambda_{\ell_1,{\ell'}_{d-1}}} d_{\underline{\ell},\ell_{d-1}} c_{\ell_1}\left(\left(\underline{\ell},\ell_{d-1}\right),\left(\underline{\ell'},{\ell'}_{d-1}\right)\right) \overline{d_{\underline{\ell'},{\ell'}_{d-1}}}\geq 0,$$
for all $ \ell_1\in \Z$, $k \geq \vert\ell_1\vert $ and arbitrary values $d_{\underline{\ell},\ell_{d-1}}\in \C$. Thereby the limit in the penultimate equation is non-negative for all $\ell_1$ and thereby the infinite sum is as well.
\end{proof}

We do not need the property of positive definiteness for arbitrary values of the index set in the proof of the result. Therefore we can loosen the assumption by introducing a matrix-like notation for the mapping. 
We first introduce a new index $\alpha$ where $\alpha\in\N$ and there is a one-to one mapping from $\N$ to the pairs 
$$\left(\underline{\ell}_{\alpha},\ell_{d-1,\alpha}\right)\in \Lambda_{\vert \ell_1\vert }$$ satisfying $\ell_{d-1,\alpha'}\leq \ell_{d-1,\alpha}$ if $\alpha' < \alpha$. 
For the proof, the actual order of the mapping is insignificant. We additionally define for each $\alpha$ the eigenvalue corresponding to $Y_{\ell_1, \underline{\ell}_{\alpha},\ell_{d-1,\alpha}}$ by $\lambda_{\alpha}$ which only depends on the value of $\ell_{d-1,\alpha}$. Further we introduce the notation $N_{\alpha}:=N_{\ell_{d-1,\alpha},d}$ for the number of spherical harmonics corresponding to the eigenvalue $\lambda_{\alpha}$.

Then we define for each $k\in \N$:

\begin{equation}\label{eq:DefMatrixAxial}
A_{\ell_1}^k=\left( c_{\ell_1}\left(\underline{\ell}_{\alpha},\ell_{d-1,\alpha},\underline{\ell}_{\alpha'},\ell_{d-1,\alpha'}\right)\right)_{\alpha,\alpha'=1}^k.
\end{equation}
 We note that in the case of $d=3$ the matrix structure of the mappings $c_{\ell_1}$ is evident since $\Lambda_{\vert\ell_1\vert}=\lbrace j\in \N \ \vert \ j\geq \vert \ell_1\vert \rbrace$.
\begin{proposition}
Let $K$ be an Hermitian, axially symmetric kernel of the form \eqref{eq:Formaxker2}, $d\geq 3$, satisfying \eqref{eq:Critabssum} for absolute summability and $A_{\ell_1}^k$ as above. 
The kernel is positive definite if the matrix $A_{\ell_1}^k$ is positive semi-definite for all $\ell_1\in\Z$ and all $k\in \N$.
\end{proposition}

The proof is similar to the proof of \Cref{Le:MapAxiaPD} and therefore omitted.

\begin{lemma}\label{Lemma:StrictposdefAxial}
For $K$ as in the above proposition.
The kernel is strictly positive definite if  $A_{\ell_1}^k$ is strictly positive definite for all $\ell_1\in \Z$ and all $k\in\N$ and if there exist sequences $d_j^{\ell_1}>0$, $j\in \N$, such that the matrix  
$$\tilde{A}_{\ell_1}^k=A_{\ell_1}^kD_k,\quad \text{with } D_k=\left(\delta_{i,j}d_j^{\ell_1}\right)_{i,j=1}^k$$  satisfies $$\lambda_{\min}\left(\tilde{A}^k_{\ell_1}\right)>\varepsilon_{\ell_1}, $$
where $\varepsilon_{\ell_1}>0$ is independent of $k$.
\end{lemma}
\begin{proof}
From the proof of the last lemma we deduce:
\begin{align*} \sum_{\xi,\zeta \in \Xi}\lambda_{\xi} \overline{\lambda_{\zeta}} K\left(\xi, \zeta\right)&=
 \sum_{\ell_1=-\infty}^{\infty}\underset{k \rightarrow \infty}{\lim} \left(\mathbf{y}^{\ell_1,k}\right)^T A^k_{\ell_1}\overline{\mathbf{y}^{\ell_1,k}}\\
 & \text{where } \mathbf{y}^{\ell_1,k}=\left(y^{\ell_1,\underline{\ell}_{\alpha},\ell_{d-1,\alpha}}_{\Xi}\right)_{\alpha=1}^{k}.
\end{align*}

Inserting the Cholesky-decomposition of $A^k_{\ell_1}=G_{\ell_1}^k\left(G_{\ell_1}^k\right)^*$ into the last equation yields
\begin{align*}
\sum_{\xi,\zeta \in \Xi }\lambda_{\xi} \overline{\lambda_{\zeta}} K\left(\xi, \zeta\right)
&=\sum_{\ell_1=-\infty}^{\infty}\underset{k \rightarrow \infty}{\lim}\left\Vert G^k_{\ell_1}\mathbf{y}^{\ell_1,k}\right\Vert_2^2\\
&=\sum_{\ell_1=-\infty}^{\infty}\underset{k \rightarrow \infty}{\lim}\left\Vert G^k_{\ell_1}\sqrt{D_k}\sqrt{D_k}^{-1}\mathbf{y}^{\ell_1,k}\right\Vert_2^2\\
&\geq \sum_{\ell_1=-\infty}^{\infty}\underset{k \rightarrow \infty}{\lim} \left(\left(\left\Vert \left(G^k_{\ell_1}\sqrt{D_k}\right)^{-1}\right\Vert_2\right)^{-1}\left\Vert\sqrt{D_k}^{-1}\mathbf{y}^{\ell_1,k}\right\Vert_2\right)^2,
\end{align*}
where we used the Cauchy-Schwarz inequality for the last estimate. 
We can now estimate the norm of the first product with the bound assumed for the eigenvalue $\tilde{A}^{k}_{\ell_1}$. We know that 
\begin{align*}\lV \left(G^{k}_{\ell_1}\sqrt{D_k}\right)^{-1}\rV_2^{-1}&=\sqrt{\lambda_{\max}\left(\sqrt{D_k}^{-1}\overline{\left(G^{k}_{\ell_1}\right)^{-1}}^T \left(G^{k}_{\ell_1}\right)^{-1}\sqrt{D_k}^{-1}\right)}^{-1}\\
&=\sqrt{\lambda_{\max}\left(\sqrt{D_k}^{-1}\left(A_{\ell_1}^{k}\right)^{-1}\sqrt{D_k}^{-1}\right)}^{-1}\\
&=\sqrt{\lambda_{\max}\left(\left(A_{\ell_1}^{k}\right)^{-1}D_k^{-1}\right)}^{-1}=\sqrt{\lambda_{\min}\left(A_{\ell_1}^{k}D_k\right)}>\sqrt{\varepsilon_{\ell_1}},
\end{align*}
independently of $k$. 
Thereby
\begin{align*}
\sum_{\xi \in \Xi} \sum_{\zeta \in \Xi}\lambda_{\xi} \overline{\lambda_{\zeta}} K\left(\xi, \zeta\right)
&\geq \sum_{\ell_1=-\infty}^{\infty}\varepsilon_{\ell_1}\underset{k \rightarrow \infty}{\lim} \lV \left(\sqrt{D_{k}}^{-1}\mathbf{y}^{\ell_1,k}\right)\rV_2^2\\
&\geq \sum_{\ell_1=-\infty}^{\infty}\varepsilon_{\ell_1} \sum_{\alpha=1}^{\infty}\left(\left(d_{\alpha}^{\ell_1}\right)^{-1}\lv y_{\Xi}^{\ell_1,\underline{\ell}_{\alpha},\ell_{d-1,\alpha}}\rv^2\right),
\end{align*}
where the last sum can only be zero if  $y_{\Xi}^{\ell_1,\underline{\ell}_{\alpha},\ell_{d-1,\alpha}}$ is the zero sequence for all $\ell_1$. This implies that the $\lambda_{\xi}$ in the definition of the quadratic form need all be zero as a consequence of the linear independence of the spherical harmonics.

\end{proof}
The conditions of the last lemma are easily verified for example for the case where $A^k_{\ell_1}$ is a diagonal matrix, this would make the kernel a convolutional kernel. In general the existence of applicable $d_j$ is not easily verified. Therefore we give a stronger condition that is easier to verify.
\begin{theorem}
Let $K$ be an Hermitian, axially symmetric kernel of the form \eqref{eq:Formaxker2} satisfying \eqref{eq:Critabssum} for absolute summability and $d>2$.
The kernel is strictly positive definite if  $$\tilde{A}^{\ell_1}=\left(\tilde{a}^{\ell_1}_{\alpha,\alpha'}\right)_{\alpha,\alpha'=1}^{\infty}=\left(\sqrt{N_{\alpha}} c_{\ell_1}\left(\underline{\ell}_{\alpha},\ell_{d-1,\alpha},\underline{\ell}_{\alpha'},\ell_{d-1,\alpha'}\right)\sqrt{N_{\alpha'}}\right)_{\alpha,\alpha'=1}^{\infty},$$ is positive definite for all $\ell_1\in \N$  and  satisfies the uniform strict diagonal dominance property:
$$\sum_{\alpha'\neq \alpha} \vert \tilde{a}^{\ell_1}_{\alpha,\alpha'}\vert < \sigma_{\ell_1} \vert \tilde{a}^{\ell_1}_{\alpha,\alpha}\vert,\quad \forall \alpha \in \N, $$
where each $0< \sigma_{\ell_1} <1$ is independent of $\alpha$. 
\end{theorem}
\begin{proof}
From the proof of the last lemma we can continue, this time we denote by $\tilde{A}^{\ell_1}=G_{\ell_1}G_{\ell_1}^*$ the infinite Cholesky decomposition matrix which exists since $\tilde{A}^{\ell_1}$ is positive definite and bounded on $\ell^{\infty}(\N)$ as a consequence of \eqref{eq:Critabssum}.

Inserting the Cholesky-decomposition into the quadratic form 
\begin{align}\label{eq:CholeskyKern}
\sum_{\xi \in \Xi} \sum_{\zeta \in \Xi}\lambda_{\xi} \overline{\lambda_{\zeta}} K\left(\xi, \zeta\right)
&=\sum_{\ell_1=-\infty}^{\infty}\lV G_{\ell_1}\left(\tilde{\mathbf{y}}^{\ell_1}\right)\rV_2^2,
\end{align}
where $\tilde{\mathbf{y}}^{\ell_1}=\left(\frac{1}{\sqrt{N_{\alpha}}}y_{\Xi}^{\ell_1,\underline{\ell}_{\alpha},\ell_{d-1,\alpha}}\right)_{\alpha=1}^{\infty}$ is a bounded sequence as a result of the estimate of the spherical harmonics in \eqref{eq:abssphH}. 
We first deduce that if there exists an element of the space of bounded sequences $\ell^{\infty}(\N)$ for which $G_{\ell_1}^*x=0$, then $G_{\ell_1}\left(G_{\ell_1}^*x\right)=\tilde{A^{\ell_1}}x=0$, where we can exchange the order of multiplication because of the triangular structure of $G_{\ell_1}$ and because $G_{\ell_1}^*x$ is elementwise finite for all bounded $x$, as a result of $\lV G^*_{\ell_1}x \rV =\overline{x}^T\tilde{A^{\ell_1}}x <\infty$ which follows from \eqref{eq:Critabssum}.

Thereby it is sufficient to show that there exists no eigenvector in the space of bounded sequences for which $\tilde{A}^{\ell_1}x=0.$
We prove by contradiction, adapting the conditions of \citep{Shivakumar1987} Theorem 1b.

From the positive definiteness of the sub-matrices of $\tilde{A}^{\ell_1}$ we can deduce that $\tilde{a}^{\ell_1}_{\alpha,\alpha}$ is non zero for all $\alpha\in\N$.

Assume there exists $x\in\ell^{\infty}(\N)$ with
$\tilde{A}^{\ell_1}x=0$ and $\lV x\rV_{\infty}=1$. Thereby 
\begin{align*} \sum_{\alpha'=0}^{\infty} \tilde{a}^{\ell_1}_{\alpha,\alpha'}x_{\alpha'}&= 0,\quad  \forall \alpha \in \N \\
\Leftrightarrow \ \sum_{\alpha'\neq \alpha}^{\infty} \tilde{a}^{\ell_1}_{\alpha,\alpha'}x_{\alpha'}&= -\tilde{a}_{\alpha,\alpha'}x_{\alpha}, \quad \forall \alpha \in \N.
\end{align*}
Since $x$ has $\ell^{\infty}$ norm one, there exists a value $\alpha\in \N$, for which $\vert x_\alpha \vert>\sigma_{\ell_1}$ and
$$\sum_{\alpha'\neq \alpha}^{\infty}\vert  \tilde{a}^{\ell_1}_{\alpha,\alpha'}\vert \geq  \sum_{\alpha'\neq \alpha}^{\infty}\vert  \tilde{a}^{\ell_1}_{\alpha,\alpha'} x_{\alpha'}\vert \geq  \vert \tilde{a}_{\alpha,\alpha}\vert \sigma_{\ell_1}.$$
This contradicts the assumption and the only possible choice for which \eqref{eq:CholeskyKern} is equal to zero is $\tilde{\textbf{y}}^{\ell_1}=0$ for all $\ell_1\in\Z$. From the definition of this sequence we know that this is only possible if $\lambda_{\xi}=0$ for all $\xi \in \Xi$ as in the proof of the last lemma.

\end{proof}
  
With the aim of finding more general sufficient conditions that are easy to evaluate we define the set $\mathcal{F}$ as the set of all indices $\ell_1 \in  \Z$ for which
$$\sum_{j,j'\in \Lambda_{\vert \ell_1\vert}}d_jc_{\ell_1}\left(j,j'\right)\overline{d_{j'}}=0\quad \Rightarrow\quad d_{j}=0\ \forall j \in \Lambda_{\vert \ell_1\vert}. $$

\begin{lemma}\label{LemmaImplSPDAxial}
A  Hermitian axially symmetric  kernel $K$ of the form \eqref{eq:Formaxker2} , satisfying \eqref{eq:Critabssum}, with positive definite maps $c_{\ell_1}:\Lambda_{\vert \ell_1\vert }\times\Lambda_{\vert \ell_1\vert}\rightarrow \C$ is strictly positive definite if
$$\sum_{\xi \in \Xi}\lambda_{\xi} Y_{\ell_1,\ell}\left(\xi\right)=0\  \forall  \ell_1\in  \mathcal{F}, \ \ell \in \Lambda_{\vert \ell_1\vert} \quad  \Rightarrow\quad \lambda_{\xi }=0\ \forall \xi \in \Xi.$$
\end{lemma}
\begin{proof}
We prove the result by contradiction and assume that $K$ is continuous positive definite but not strictly positive definite.  If $K$ is not strictly positive definite there exists a nonempty set of distinct point $\Xi$ and coefficients $\lambda_{\xi}$ not all zero with  
$$ \sum_{\xi,\zeta \in \Xi}\lambda_{\xi} \overline{\lambda_{\zeta}} K\left(\xi, \zeta\right)=0. $$
This is equivalent according to the computation in the proof of \Cref{Le:MapAxiaPD} to 
\begin{equation}\label{eq:Formaxker3}
\sum_{\ell_1=-\infty}^{\infty}\sum_{\ell,\ell' \in \Lambda_{\vert \ell_1\vert }}c_{\ell_1}\left(\ell,\ell'\right) Y_{\ell_1,\ell,\Xi} \overline{Y_{\ell_1,\ell',\Xi}}=0,
\end{equation}
where $Y_{\ell_1,\ell,\Xi}=\sum_{\xi\in\Xi} \lambda_{\xi} Y_{\ell_1,\ell}\left(\xi\right)$.

Since we know that all the sums are non negative since the maps  $c_{\ell_1}$ are all positive definite, the overall sum can only be zero if all summands are. For the indices $\ell_1\in \mathcal{F}$ this implies  $Y_{\ell_1,\ell,\Xi}=0$.
But according to the condition of the lemma this implies $\lambda_\xi=0$ for all $\xi\in\Xi$, in contradiction to the assumption of the proof.
\end{proof}
The following proposition is an immediate consequence of the last lemma together with the linear independence of the spherical harmonics.
\begin{proposition}
A kernel as in in the last lemma is strictly positive definite if $\mathcal{F}=\Z$.
\end{proposition}
\begin{lemma}\label{Lemma:SPDAxialNec}
Let $K$ be an Hermitian, axially symmetric kernel of the form \eqref{eq:Formaxker2}, satisfying \eqref{eq:Critabssum}.
For the kernel to be strictly positive definite it is necessary that the mapping
 $$c_{\ell_1}: \Lambda_{\vert \ell_1\vert }\times \Lambda_{\vert \ell_1 \vert } \rightarrow \C$$
 is not identically zero for infinitely many $\ell_1 \in \Z$. 
\end{lemma}
\begin{proof}
The result  is proven by contradiction. Assume there are only finitely many values of $\ell_1$ for which $c_{\ell_1}$ is not the zero mapping.
Let us denote the set of these $\ell_1$ by $\mathcal{J}$. Then we can construct a set of points $\Xi$ and coefficients $\lambda_{\xi}\in \R$ such that
\begin{equation}
\sum_{\xi\in \Xi}\lambda_{\xi}Y_{\ell_1,\ell}\left(\xi\right)=0,\quad \forall \ell_1\in \mathcal{J}, \ \ell \in \Lambda_{\vert \ell_1\vert }.
\end{equation}
These are infinitely many conditions since $ \Lambda_{\vert \ell_1\vert }$ includes infinitely many elements but if we choose all points in $\Xi$ of the form $\left(\theta_{1,k}, \theta_{2}, \ldots,\theta_{d-1}\right)$  for $k=0,\ldots,\vert \mathcal{J}\vert$, the above equations can be transformed into
\begin{equation}
\frac{1}{\sqrt{2\pi}}\prod_{j=2}^{d-1}{}_{j}\tilde{P}_{\ell_j}^{\ell_{j-1}}\left(\theta_{n}\right)\sum_{k=0}^{\vert \mathcal{J}\vert}\lambda_{k}e^{i\ell_1\theta_{1,k}} =0,\quad \forall \ell_1\in\mathcal{J}, \ \ell\in \Lambda_{\vert \ell_1\vert},
\end{equation}
which is satisfied if 
\begin{equation}
\sum_{k=0}^{\vert\mathcal{J}\vert}\lambda_{k}e^{i\ell_1\theta_{1,k}} =0,\quad \forall \ell_1\in\mathcal{J}.
\end{equation}
These are only  $\vert\mathcal{J}\vert$ conditions which can be satisfied for our choice of $\theta_{k,1}$ at  $\vert\mathcal{J} \vert+1$ distinct points.
\end{proof}


\begin{lemma}\label{Lemma:SPDAxialNec2}
Let $K$ be an Hermitian, axially symmetric kernel of the form \eqref{eq:Formaxker2} satisfying \eqref{eq:Critabssum}.
For the kernel to be strictly positive definite, it is necessary that the mapping
 $$c_{0}: \Lambda_{\vert \ell_1 \vert }\times \Lambda_{\vert \ell_1\vert } \rightarrow \C$$
 is not constant zero. 
\end{lemma}
\begin{proof}
We show that for any point on the surface of the sphere which has the form $\left(\theta_1,0,\theta_3,\ldots,\theta_{d-1}\right)=\xi$, the spherical basis functions $Y_{\ell_1,\ell}\left(\xi\right)$ only take non zero values if $\ell_1=0$. Inserting the point into the definition of the spherical harmonics from \eqref{eq:spherharmaxial}
\begin{align*}
Y_{\ell_1,\ell}\left(\xi\right)&=\frac{1}{\sqrt{2\pi}}e^{i\ell_1\theta_1} \prod_{j=2}^{d-1}{}_{j}\tilde{P}_{\ell_j}^{\ell_{j-1}}\left(\theta_{n}\right)\\
&=\frac{1}{\sqrt{2\pi}}e^{i\ell_1\theta_1}\left( \frac{2\ell_2+1}{2}\frac{\left(\ell_2+\ell_1\right)!}{\left(\ell_2-\ell_1\right)!}\right)^{1/2}P_{\ell_2}^{-\ell_1}\left(1\right) \prod_{j=3}^{d-1}{}_{j}\tilde{P}_{\ell_j}^{\ell_{j-1}}\left(\theta_{j}\right).
\end{align*}
Since the last part of the product is finite, the value is non-zero only if $P_{\ell_1}^{-\ell_2}\left(1\right)$ is non-zero and this is only the case for $\ell_1=0$.
\end{proof}
The study of convolutional kernels in the next sections will give additional results on strict positive definiteness of axially symmetric kernels especially for the case where $d=2$.
\section{(Strict) positive definiteness of convolutional kernels or kernels with eigenvalue block structure}
We now assume that $K$ is continuous Hermitian and has the form  
\begin{equation}\label{eq:Kernform2} K\left(\xi, \zeta\right)=\sum_{j=0}^{\infty}\sum_{k=1}^{N_{j,d}} \sum_{k'=1}^{N_{j,d}} d_j\left(k,k'\right) Y_j^k\left(\xi\right)\overline{Y_{j}^{k'}\left(\zeta\right)},\qquad \forall \xi,\zeta \in \S,
 \end{equation}
 for a fixed set of orthogonal basis functions, where $Y_{j}^{k}$ is an eigenfunction  corresponding to the eigenvalue $\lambda_j$ of \eqref{eq:EigenProb}. We note that under this assumption there exists another orthonormal basis of the eigenfunctions corresponding to each eigenvalue such that, if we express $K$ with respect to this sequence (which is also the Hilbert-Schmidt basis of the kernel), it takes the form \eqref{eq:Kernform3}. 
 
These kernels were studied under the name convolutional kernels for example in \citep{Narcowich1995},\citep{Narcowich2002}. The fixed basis case on the 2-sphere was also studied in \citep{Jaeger2021} and we now generalise the results.

It is proven in \citep{Dyn1999} that a continuous and Hermitian kernel $K$ of the form \eqref{eq:Kernform3} is positive definite if and only if the $d_{j,k}$ are real and non-negative for all  $j \in \Z_+$, $k\in\left\lbrace 1,\ldots,N_{j,d}\right\rbrace$, ( their Theorem 2.1).

We can now transfer the result known for radial kernels which was first proven in \citep{Chen2003}, to this kernel class. Since positive definiteness is independent of the choice of the expansion basis, the proofs will be carried out using kernels of the form \eqref{eq:Kernform3} and we later briefly note the implications for fixed basis kernels.

\begin{lemma}\label{LemmaInfEvenInfOdd} If a continuous and Hermitian kernel $K$  of the form \eqref{eq:Kernform3}
 is strictly positive definite, all $d_{j,k}$ are non-negative and the sequences $\left(d_{j,k}\right)_{k=1}^{N_{j,d}}$ are not identically zero for infinitely many even and infinitely many odd values of $j\in \Z_+$.
\end{lemma}
\begin{proof}
We assume $K$ is continuous, Hermitian and strictly positive definite. Since strict positive definiteness implies positive definiteness, the $d_{j,k}$ are all non negative according to \citep{Dyn1999}, Theorem 2.  
The rest of the proof is divided into four cases. Let us denote the set of indices for which $\left(d_{j,k}\right)_{k=1}^{N_{j,d}}$ are  not all zero with $$\mathcal{J}:=\left \lbrace j \in \Z_+ \vert \left(d_{j,1},\ldots,d_{j,N_{j,d}}\right)^T \neq \mathbf{0}^{N_{j,d} }\right \rbrace.$$
For all fourcases  we assume that $K$ is strictly positive definite of the form  \eqref{eq:Kernform3} and prove that assuming  either
\begin{enumerate}
\item $\mathcal{J} \subset 2\N$ or
\item $\mathcal{J}\subset 2\N+1$ or
\item $1\leq \vert \mathcal{J}\cap 2\N\vert< \infty$ or
\item $1\leq \vert \mathcal{J}\cap \left(2\N+1\right)\vert< \infty$,
\end{enumerate}
leads to a contradiction.

1. Let us now assume $\left(d_{j,k}\right)_{k=1}^{N_{j,d}}$ are  not all zero \textbf{ only for even} values of $j$. We can represent the quadratic form as 
\begin{align*} \sum_{\xi,\zeta \in \Xi}\lambda_{\xi} \overline{\lambda_{\zeta}} K\left(\xi, \zeta\right)&=\sum_{j=0}^{\infty}\sum_{k=1}^{N_{j,d}}d_{j,k}\sum_{\xi \in \Xi}\lambda_{\xi} \tilde{Y}_j^k\left(\xi\right) \overline{\sum_{\zeta \in \Xi}\lambda_{\zeta}\tilde{Y}_{j}^{k}\left(\zeta\right)} \\
& =\sum_{j=0}^{\infty}\sum_{k=1}^{N_{j,d}} y^{j,k}_{\Xi} d_{j,k}\overline{y^{j,k}_{\Xi}},\quad \text{ with } y^{j,k}_{\Xi}=\sum_{\xi \in \Xi}\lambda_{\xi} \tilde{Y}_j^k\left(\xi\right) \in\C.
\end{align*}
Choosing a non-empty set of data sites $\Xi'$ which satisfies 
$$\xi \in \Xi'\ \Rightarrow \ -\xi \in \Xi',$$
where $-\xi$ is the antipodal point of $\xi$, and setting $\lambda_{\xi}=-c_{-\xi}$, we find that $y_{\Xi'}^{j,k}=0$ for all even $j$ since $\tilde{Y}_j^k\left(-\xi\right)=\left(-1\right)^j\tilde{Y}_j^k\left(\xi\right)$ follows from the $\tilde{Y}_j^k$ all being homogeneous polynomials of order $j$. 
This implies $\sum_{\xi,\zeta \in \Xi'}\lambda_{\xi} \overline{\lambda_{\zeta}} K\left(\xi, \zeta\right)=0$ and therefore $K$ would not be strictly positive definite.

2. The same argument applies when we assume $\left(d_{j,k}\right)_{k=1}^{N_{j,d}}$ are  not all zero \textbf{ only for odd} values of $j$. The same set $\Xi'$ can be chosen but  now $\lambda_{\xi}=\lambda_{-\xi}$ yields the contradiction to the assumption of strict positive definiteness.

3. Now we assume that $\left(d_{j,k}\right)_{k=1}^{N_{j,d}}$ are  not all zero for any number of odd values of $j$ and\textbf{ only finitely many even} values of $j$. Set $\hat \jmath$ to the maximal even index for which $\left(d_{j,k}\right)_{k=1}^{N_{j,d}}$ is not zero.
We aim to construct a set $\Xi$ with elements only in the lower hemisphere of $\S$ and $\lambda\in \C^{\vert \Xi \vert}\neq 0$, s.t. 
$\sum_{\xi \in \Xi}\lambda_{\xi} \tilde{Y}_j^k\left(\xi\right) =0$ for all even $j \leq \hat\jmath$. 
These are $$M:=\sum_{m=1}^{\hat\jmath/2}N_{2m,d}<\infty$$
 linear equations.
We can therefore choose any set of distinct points in the lower hemisphere with more than $M$ elements to find a non trivial solution.  Defining the set $\Xi'=\Xi \cup\left( -\Xi\right)$ and setting $\lambda_{\xi}=-\lambda_{-\xi}$ shows that 
 $$\sum_{\xi,\zeta \in \Xi'}\lambda_{\xi} \overline{\lambda_{\zeta}} K\left(\xi, \zeta\right)=0$$
 for a non trivial vector $\lambda$ and therefore $K$ is not strictly positive definite.
 
4. The same arguments can be used to show that $\left(d_{j,k}\right)_{k=1}^{N_{j,d}}$ needs to be\textbf{ nonzero for infinitely many odd} values of $j$.

\end{proof}

For the fixed basis form this implies that for a strictly positive definite kernel of the form \eqref{eq:Kernform2} it is necessary that the matrices $D_j:=\left(d_j\left(k,k'\right)\right)_{k,k'=1}^{N_{j,d}}$ are all positive definite and infinitely many with odd $j$ and infinitely many with even $j$ are not the zero matrix.

With the aim of finding sufficient conditions that are easy to evaluate we define the set $\mathcal{F}$ as the set of all indices $j \in  \Z_+$ for which $D_j$ is a  positive definite matrix or if the appropriate basis is used all $d_{j,k}>0$ for all $k=1,\ldots,N_{j,d}$.
\begin{lemma}\label{LemmaImplSPD}
A continuous and Hermitian kernel $K$ of the form \eqref{eq:Kernform2} which is positive definite is strictly positive definite if
$$\sum_{\xi \in \Xi}\lambda_{\xi} Y_j^k\left(\xi\right)=0,\  \forall \ k=1, \ldots,N_{j,d},\  j\in \mathcal{F} \quad  \Rightarrow\quad \lambda_{\xi }=0, \forall \xi \in \Xi.$$
\end{lemma}
The proof follows in the same line of argument as \Cref{LemmaImplSPDAxial}.


With \Cref{LemmaImplSPD}, we have shown that strict positive definiteness of these kernels can be proven using the same results which were used for zonal kernels by Chen et al.\ in \citep{Chen2003} and the alternative proof for more general manifolds stated by Barbosa and Menegatto in \citep{Barbosa2016}. The next lemma and \Cref{THM:SphereSuf}, generalise the results to the convolutional case and complex kernels.
\begin{lemma}\label{Lemma:EquivCondSPD}
For a given set $\mathcal{F}$ the following two properties are equivalent:
\begin{enumerate}
\item $\sum_{\xi \in \Xi}\lambda_{\xi} Y_j^k\left(\xi\right)=0,\  \forall \ k=1, \ldots,N_{j,d},\  j\in \mathcal{F}$ implies  $\lambda_{\xi }=0,\ \forall \xi \in \Xi.$
\item $ \sum_{\xi\in \Xi }\lambda_{\xi} P^d_{j}\left(\xi^T\zeta\right)=0,\ \forall j \in \mathcal{F},\, \forall \zeta \in \S$ 
implies $\lambda_{\xi }=0,\ \forall \xi \in \Xi$.
\end{enumerate}Here $P_j^{d}$ are the Legendre polynomials of degree $j$ in $d$ dimensions.
\end{lemma}
\begin{proof}
Using the addition formula of the Legendre polynomials we find $$\frac{N_{j,d}}{\sigma_{d-1}}P_{j}^d\left(\xi^T\zeta\right)=\sum_{k=1}^{N_{j,d}} Y_{j}^k\left(\xi\right)\overline{Y_{j}^k\left(\zeta\right).}$$
Therefore 
$$\frac{N_{j,d}}{\sigma_{d-1}} \sum_{\xi \in \Xi}\lambda_{\xi} P_{j}^d\left(\xi^T\zeta\right)=\sum_{k=1}^{N_{j,d}} \left(\sum_{\xi \in \Xi}\lambda_{\xi}Y_{j}^k\left(\xi\right)\right) \overline{Y_{j}^k\left(\zeta\right)}, \qquad \forall \zeta\in \S.$$
Since the spherical harmonics are linearly independent, the last expression is  zero if and only if 
$$\sum_{\xi \in \Xi}\lambda_{\xi}Y_{j}^k\left(\xi\right)=0$$ for all $k=1,\ldots ,N_{j,d}$, and all $j\in \mathcal{F}$. 
\end{proof}

\begin{theorem}\label{THM:SphereSuf}
Let $K$ be a continuous positive definite kernel of the form \eqref{eq:Kernform2}, $d>2$ and $\mathcal{F}$ the corresponding index set for which the $D_j$ are positive definite matrices. Then it is sufficient for $K$ to be strictly positive definite that $\mathcal{F}$ includes infinitely many even and infinitely many odd values of $j\in \Z_+$.
\end{theorem} 
\begin{proof}Combining \Cref{LemmaInfEvenInfOdd} and \Cref{LemmaImplSPD} we know that for $K$ to be strictly positive definite  it is sufficient to prove that for arbitrary sets of distinct data sites $\Xi \subset \S$ the functions 
$P_j^d\left(\xi^T\zeta\right)$ satisfy \Cref{LemmaImplSPD} (2). Therefore we show that for any such set $\Xi$
$$\sum_{\xi \in \Xi}\lambda_{\xi}P_j^d\left(\xi^T\zeta\right)=0, \qquad  \forall j \in \mathcal{F},\ \zeta\in \S,$$
implies $\lambda_{\xi}=0$ for all  $\xi \in \Xi$.
We do this by choosing for each $\xi\in \Xi$ a corresponding $\zeta=\zeta_{\xi}\in \S$ and show that for this choice $\lambda_{\xi}=0$ if the above holds.
Assume we have fixed $\xi \in \Xi$, we need to distinguish two cases.

\textbf{Case 1:} $ \xi^T\zeta\neq -1$ for all $\zeta \in \Xi.$

In this case we choose $\zeta_{\xi}=\xi$ and the system above takes the form  
$$\lambda_{\xi} P_j^d\left(1\right) +\sum_{\zeta \in \Xi\setminus \lbrace \xi \rbrace}\lambda_{\zeta}P_j^d\left(\zeta^T\xi\right)=0, \qquad  \forall j \in \mathcal{F}.$$
Note that $\xi^T\zeta \in \left(-1,1\right)$ and we set $\theta_{\zeta}=\arccos\left(\xi^T\zeta\right)\in \left(0,\pi\right)$. 
Since there are infinitely many indices in $\mathcal{F}$ we can choose a sequence from $\mathcal{F}$ with $j_n\in \mathcal{F}$ and $\underset{n \rightarrow \infty}{\lim} j_n=\infty$.
Introducing the limit in the above equation and using the asymptotic estimate for Legendre polynomials in \citep{Atkinson2012} (2.117) 
\begin{align*} \lambda_{\xi} + \underset{n \rightarrow \infty}{\lim} \frac{\sum_{\zeta \in \Xi\setminus \lbrace \xi \rbrace}\lambda_{\zeta}P_j^d\left(\theta_{\zeta}\right)}{P_j^{d}\left(1\right)}=0 
\end{align*}
implies $\lambda_{\xi}=0$ because the sum is finite and the latter part converges to zero.

\textbf{Case 2} There is one $\zeta \in \Xi $ with $\xi^T \zeta =-1$.

Then $\zeta$ is the antipodal point of $\xi$ and in Cartesian coordinates $\xi =-\zeta$.  Then the above equation becomes 
$$\lambda_{\xi} P_j^d\left(1\right) +\left(-1\right)^j \lambda_{-\xi} P_j^d\left(1\right)+\sum_{\zeta \in \Xi\setminus \lbrace \xi, -\xi \rbrace}\lambda_{\zeta}P_j^d\left(\zeta^T\xi\right)=0, \qquad  \forall j \in \mathcal{F}.$$
The third part vanishes if we introduce the limit as in the previous case but for odd and even series of $j_n$ separately. The remainder yields for even $j$ and odd $j$, respectively,
$$ \lambda_{\xi }+\lambda_{-\xi}=0= \lambda_{\xi}-\lambda_{-\xi}$$
which implies $\lambda_{\xi}=0=\lambda_{-\xi}$.
\end{proof}

\section{The special case of axially symmetric kernels on the circle}
For the remaining case $d=2$ we note that a basis of the spherical harmonics is
\[Y_{j}\left(\theta\right)=\frac{1}{\sqrt{2\pi}}e^{ij\cdot\theta},\quad j\in \Z,\]
where $\theta$ is associated with the point on the circle $\xi=\left(\cos\left(\theta\right),\sin\left(\theta\right)\right)^T$. The Legendre polynomials take the form of Chebyshev polynomials
\[P_j^2\left(\xi^T\zeta\right)=\cos\left(j \xi^T\zeta\right).\]
Thereby every continuous kernel on $\mathbb{S}^1\times \mathbb{S}^1$  has an expansion of the form
\[ K(\xi,\zeta)=\sum_{j,j'\in\Z}a_{j,j'}e^{ij\cdot \theta}e^{-ij'\cdot\theta'},\]
where $\theta$ is corresponding to $\xi$ as above and $\zeta=\left(\cos\left(\theta'\right),\sin\left(\theta'\right)\right)^T$.
According to \Cref{THM:CoefStrucAX}, such a kernel is axially symmetric if $$a_{j,j'}=\delta_{j,j'}c_j.$$
For an Hermitian convolutional kernel only the coefficients corresponding to the same order of trigonometric polynomial are going to be non-zero. This yields 
\begin{align*} K(\xi,\zeta)=&\sum_{j\in \Z_+}\left(d_{j,j}e^{ij(\theta-\theta')}+d_{-j,-j}e^{-ij(\theta-\theta')}+d_{-j,j}e^{-ij(\theta+\theta')}+\overline{d_{-j,j}}e^{ij(\theta+\theta')}\right)\\ &+d_{0,0},\end{align*}
where we can define the matrices of the last section as $D_j=\begin{pmatrix}
d_{j,j} & d_{j,-j} \\ \overline{d_{j,-j}} &d_{-j,-j}
\end{pmatrix}$.
\begin{theorem}\label{THM:CircelSuf}
Let $K$ be a continuous positive definite kernel of the form \eqref{eq:Kernform2}, $d=2$, and $\mathcal{F}$ the corresponding index set for which the $D_j$ are positive definite matrices. Then it is sufficient for $K$ to be strictly positive definite that $\mathcal{F}$ intersects non-trivially every full arithmetic progression in $\Z_+$.
\end{theorem}
\begin{proof}
For this special case \Cref{Lemma:EquivCondSPD} (2) reads
\[\sum_{\xi\in\Xi}\lambda_{\xi}\cos\left(j\arccos\left(\xi^T\zeta\right)\right)=0,\quad \forall j \in \mathcal{F},\ \forall \xi \in \mathbb{S}^1\]
implies $\lambda_{\xi}=0$ for all $\xi \in \Xi$.
This is the case if all functions of the form $\sum_{j\in \mathcal{F}}a_j\cos\left(j\arccos\left(\xi^T\zeta\right)\right)$ with $a_j>0$  are positive definite and thereby we can deduce our result directly from \citep{Menegatto1995}, Theorem 2.1.
\end{proof}

\begin{theorem}\label{THM:CircelNec}
Let $K$ be a continuous positive definite kernel of the form \eqref{eq:Kernform2}, $d=2$ and $\mathcal{J}$ the corresponding index set for which the $D_j$ are non zero. Then it is necessary for $K$ to be strictly positive definite that $\mathcal{J}$ intersects every full arithmetic progression in $\Z_+$.
\end{theorem}
\begin{proof}
The construction of sets which show the necessity of this condition can again be deduced from the discussion of such sets for radial kernels in \citep{Sun2005}. For radial kernels such sets were constructed for which  \Cref{LemmaImplSPD} (2) does not hold, see \citep{Barbosa2016} Theorem 5.2, and this also implies that such kernels are not strictly positive definite.
\end{proof}
\section{Discussion}

The discussion in Section 2 allows to construct kernels with different geometric properties. Even though these kernels could only be constructed approximately as truncated series we believe that they can nevertheless be used in applications as discussed in \citep{Alegria2020} for axially symmetric kernels. We believe that taking into account a wider range of geometric properties of the kernel will allow the improvement of approximation results for many areas of application. Unfortunately any truncation would remove strict positive definiteness.

The more specific results in Section 3 show that stating necessary and sufficient conditions for axially symmetric strictly positive definite kernels is possible in a similar way to the results known for radial kernels and that in fact one can utilise the results known for isotropic kernels to proof the more general  case, as done in Section 4.
 

\subsection*{Acknowledgement}
The work of J. J\"ager was supported by the Justus Liebig University's postdoctoral fellowship Just'Us.

\def\ln{\log}

\bibliographystyle{abbrvnat}
\bibliography{/home/janin/JLUbox/Dokumente/LiteraturALL}

\end{document}